\documentclass[11pt]{amsart}

\usepackage{amsmath,amssymb,amsfonts,graphics}
\newcommand{\A}{\mathfrak{A}}
\newcommand{\B}{\mathfrak{B}}

\newcommand{\N}{\mathbb{N}}

\newcommand{\C}{\mathbb{C}}

\newcommand{\dss}{\displaystyle}
\newcommand{\from}{\colon}
\newcommand{\supp}{\text{supp}}
\newcommand{\hT}{\text{hyper-Tauberian}}

\newcommand{\hl}{\text{hyperlocal}}
\newcommand{\hla}{\text{hyperlocal}}

\newcommand{\la}{\langle}
\newcommand{\ra}{\rangle}

\newcommand{\tmod}{\text{mod}}

\newtheorem{thm}{Theorem}[section]
\newtheorem{cor}[thm]{Corollary}

\newtheorem{prop}[thm]{Proposition}
\theoremstyle{definition}
\newtheorem{defn}[thm]{Definition}
\theoremstyle{remark}

\newtheorem{rem}[thm]{Remark}
\numberwithin{equation}{section}

\title[On local properties of Hochschild cohomology of a C$^*$- algebra]
{On local properties of Hochschild cohomology of a C$^*$- algebra}

\author[Ebrahim Samei]{Ebrahim Samei}
\address{Ebrahim Samei, EPFL-SB-IACS, Station 8, Ch-1015 Lausanne, 
Suisse}


\thanks{}

\email{ebrahim.samei@epfl.ch}

\begin{document}

\maketitle

\begin{abstract}
Let $A$ be a C$^*$-algebra, and let $X$ be a Banach $A$-bimodule.
B. E. Johnson showed that local derivations from $A$ into $X$ are derivations. 
We extend this concept of locality to the higher cohomology of a $C^*$-algebra
%for $n$-cocycles from  $A^{(n)}$ into $X$ 
and show that, for every $n\in \ \N$, bounded local  $n$-cocycles from  $A^{(n)}$ 
into $X$ are $n$-cocycles. 
\end{abstract}

The study of the local properties of Hochschild cohomology of a Banach algebra
was initiated by introducing the concept of ``local derivations". Let
$A$ be a Banach algebra, and let $X$ be a Banach $A$-bimodule. An
operator $D \from A \to X$ is a local derivation if for each $a\in
A$, there is a derivations $D_a \from A \to X$ such that
$D(a)=D_a(a)$. This concept was introduced independently 
by R. V. Kadison \cite{K} and D. R. Larson \cite{L} and it has been interests
of studies since then. Kadison's motivation was based  on his and J. R. Ringrose's earlier
investigation of Hochschild cohomology of various operator algebras,
whereas Larson's motivation was to investigate algebraic reflexivity (resp. reflexivity) of the
linear space of derivations (resp. bounded derivations) from a Banach algebra.
Local derivations have been investigated for various classes of Banach algebras
such as operator algebras, Banach operator algebras, group
algebras, and Fourier algebras (see \cite{HL}, \cite{S2}, \cite{S3} and the references therein).

In \cite{K}, Kadison showed that bounded local derivations 
from a von Neumann algebra into any of its dual bimodules are derivations. 
He then raised the question of whether the preceding 
result can be extended to the local higher cohomology. The purpose of this article 
is to answer affirmatively to this question in more general setting.  We 
show that if $A$ is a $C^*$-algebra and $n\in \N$, then bounded {\it local $n$-cocylces}
from  $A^{(n)}$ into any Banach $A$-bimodule are $n$-cocycles. This has already been obtained by B. E. Johnson in \cite{J} for the case $n=1$. Our approach 
is as follow:

Let $A$ be a Banach algebra, let $X$ be a Banach $A$-bimodule, and let $n\in \N$. 
In Section \ref{S:loc n-cocycle-Def}, we introduce certain $n$-linear maps from $A^{(n)}$ into $X$
which are more general than local $n$-cocycles. We call them {\it $n$-hyperlocal maps}.
We show that in order to characterize bounded local 
$n$-cocycles from $A^{(n)}$ into $X$, it suffices to first extend them to
$A^{\sharp (n)}$, where $A^\sharp$ is the unitization of $A$, and
view them as $n$-hyperlocal maps from $A^{\sharp (n)}$ into $X$.
Then, by imposing certain conditions, one can obtain the result. As it is shown
in Proposition \ref{P:App. Loc. n cocycle}, the advantage of this technique 
is that we can ``transfer the information" from the lower cohomology to the higher one if we
consider $n$-hyperlocal maps rather than local $n$-cocycles.

In Section \ref{S:HL}, we apply these ideas to hyper-Tauberian algebras. 
These algebras were introduced and studied in
\cite{S3} because of their useful local properties. 
By using the results of the preceding section, together
with the properties of \hT\ algebras, we show that bounded 
local n-cocycles from $A^{(n)}$ into $X$ are n-cocycles when $A$ is 
a \hT\ algebra. 

In Section \ref{S:C alg}, we first show that every commutative $C^*$-algebra
is a hyper-Tauberian algebra. We then apply the results of Section \ref{S:HL} to obtain our 
result for a general C$^*$-algebra. Finally, in the last section, we give a 
characterization of amenable C$^*$-algebras in terms of the 1-hyperlocal 
maps.

\section{Preliminaries}\label{S:Prelimin}

Let $X$ and $Y$ be Banach spaces. For $n\in \N$, let $X^{(n)}$ be
the Cartesian product of $n$ copies of $X$, and let $L^n(X,Y)$ and
$B^n(X,Y)$ be the spaces of n-linear maps and bounded n-linear
maps from $X^{(n)}$ into $Y$, respectively.

Let $A$ be a Banach algebra, and let $X$ be a Banach $A$-bimodule. An operator $D\in L(A,X)$ is a
derivation if for all $a,b\in A$, $D(ab)=aD(b)+D(a)b$. For each
$x\in X$, the operator $ad_x\in B(A,X)$ defined by $ad_x(a)=ax-xa$
is a bounded derivation, called an inner derivation. Let $Z^1(A,X)$
and $\mathcal{Z}^1(A,X)$ be the linear spaces of derivations and
bounded derivations from $A$ into $X$, respectively. For $n\in \N$
and $T\in L^n(A,X)$, define
\begin{eqnarray*}
\delta^nT:(a_1,\ldots ,a_{n+1}) &\mapsto & a_1T(a_2,\ldots ,a_n)\\
&+ & \sum_{j=1}^n (-1)^jT(a_1,\ldots ,a_{j-1},a_ja_{j+1},\ldots ,a_{n+1})\\
& + &  (-1)^{n+1}T(a_1,\ldots ,a_n)a_{n+1}.
\end{eqnarray*}
It is clear that $\delta^n$ is a linear map from $L^n(A,X)$ into
$L^{n+1}(A,X)$; these maps are the {\it connecting maps}. The
elements of $\ker \delta^n$ are the {\it n-cocycles}; we denote this
linear space by $Z^n(A,X)$. If we replace $L^n(A,X)$ with $B^n(A,X)$
in the above, we will have the `Banach' version of the connecting
maps; we denote them with the same notation $\delta^n$. In this
case, $\delta^n$ is a bounded linear map from $B^n(A,X)$ into
$B^{n+1}(A,X)$; these maps are the {\it bounded connecting maps}.
The elements of $\ker \delta^n$ are the {\it bounded n-cocycles}; we
denote this linear space by $\mathcal{Z}^n(A,X)$. It is easy to
check that $Z^1(A,X)$ and $\mathcal{Z}^1(A,X)$ coincide with our
previous definition of these notations.

Let $A$ be a Banach algebra, and let $X$ be a Banach $A$-bimodule.
By \cite[Section 2.8]{D}, for $n\in \N$, the Banach space $B^n(A,X)$
turns into a Banach $A$-bimodule by the actions defined by:
\begin{eqnarray*}
(a\star T)(a_1,\ldots , a_n) &= & aT(a_1,\ldots , a_n);
\end{eqnarray*}
\begin{eqnarray*}
(T\star a)(a_1,\ldots ,a_n) &= & T(aa_1,\ldots ,a_n)  \\
&+ & \sum_{j=1}^n (-1)^jT(a,a_1,\ldots ,a_ja_{j+1},\ldots ,a_n)\\
& + & (-1)^{n+1}T(a,a_1,\ldots ,a_{n-1})a_n.
\end{eqnarray*}
In particular, when $n=1$, $B(A,X)$ becomes a Banach $A$-bimodule
with respect to the products
$$(a\star T)(b)=aT(b) \ \ , \ (T\star a)(b)=T(ab)-T(a)b.$$
Let $\Lambda_n \from B^{n+1}(A,X) \to B^n(A,B(A,X))$ be the
identification given by
$$(\Lambda_n(T)(a_1,\ldots ,a_n))(a_{n+1})=T(a_1,\ldots
,a_{n+1}).$$ Then $\Lambda_n$ is an $A$-bimodule isometric
isomorphism. If we denote the connecting maps for the complex
$B^n(A,(B(A,X),\star))$ by $\Delta^n$, then we can show that
$$\Lambda_{n+1}\circ \delta^{n+1}=\Delta^n\circ \Lambda_n.$$

\section{n-hyperlocal maps and local n-cocycles}\label{S:loc n-cocycle-Def}

\begin{defn}
Let $A$ be a Banach algebra, and let $X$ be a Banach $A$-bimodule.
For $n\in \N$, let $T$ be an $n$-linear map from $A^{(n)}$ into $X$.\\
$($i$)$ $T$ is n-{\it hyperlocal} if, for $a_0,\ldots, a_{n+1}\in A$,
$$a_0a_1=a_1a_2=\cdots=a_na_{n+1}=0 \ \ \text{implies}\ \ a_0T(a_1,\ldots, a_n)a_{n+1}=0.$$
For $n=1$, 1-\hl\ maps are simply called \hl\ maps or \hl\ operators.\\
$($ii$)$ $T$ is a {\it local n-cocycle} if, for each
$\tilde{a}=(a_1,\ldots , a_n)\in A^{(n)}$, there is 
an $n$-cocycle $T_{\tilde{a}}$ from $A^{(n)}$ into $X$ such that
$T(\tilde{a})=T_{\tilde{a}}(\tilde{a})$. If, in addition, $T$ is bounded, we
say that $T$ is a {\it bounded local n-cocycle}.
\end{defn}

It is easy to see that every (local) $n$-cocycle is a $n$-hyperlocal map.
The following proposition states some sufficient conditions for a
bounded $n$-hyperlocal map to be an $n$-cocycle. This is critical for us
to obtain our result. 

\begin{prop}\label{P:App. Loc. n cocycle}
Let $A$ be a unital Banach algebra with unit 1 which
satisfies the following two conditions:\\
$($i$)$ For every unital Banach $A$-bimodule $X$, a bounded operator
$D \from A \to X$ is a left multiplier if
and only if $ba=0$ implies $D(b)a=0$.\\
$($ii$)$ For every unital Banach $A$-bimodule $X$, a bounded
operator $D \from A \to X$ is \hl\ if and only if
$$D(acb)-aD(cb)-D(ac)b+aD(c)b=0$$ for all $a,b,c \in A$.\\
Let $X$ be a unital Banach $A$-bimodule, let $n\in \N$, and let
$T\in B^n(A,X)$ be an $n$-\hl\ map such that $T(a_1,\ldots, a_n)=0$ if
any one of $a_1,\ldots, a_n$ is 1. Then $T$ is an $n$-cocycle.
\end{prop}

\begin{proof}
We prove the statement by induction on $n$. For $n=1$, by hypothesis,
$$T(acb)-aT(cb)-T(ac)b+aT(c)b=0$$ for all $a,b,c \in A$.
Since $T(1)=0$, by putting $c=1$ we get the result. Now suppose
that the result is true for $n=k$ ($k\geq 1$). We show that it is
also true for $n=k+1$. Let $T\in B^{k+1}(A,X)$ be $k+1$-\hl\ such
that $T(a_1,\ldots, a_{k+1})=0$ if any one of $a_1,\ldots,
a_{k+1}$ is 1. We first show that $\Lambda_k(T)\in B^k(A,B(A,X))$
is $k$-\hl. Let $a_0,\ldots, a_{k+1}\in A$ such that
$a_0a_1=\cdots=a_ka_{k+1}=0$, and put
$$S=a_0\star \Lambda_k(T)(a_1,\ldots, a_k)\star a_{k+1}.$$ Then $S
\from A \to X$ is a bounded operator. We claim that $S$ satisfies
the following condition:
$$bc=0 \ \text{implies}\ S(b)c=0. \eqno{(1)} $$
Let $b,c\in A$ such that $bc=0$. Then
\begin{eqnarray*}
S(b)c &=& [a_0\star \Lambda_k(T)(a_1,\ldots, a_k)\star
a_{k+1}](b)c
\\ &= & a_0(\Lambda_k(T)(a_1,\ldots, a_k))(a_{k+1}b)c
- a_0(\Lambda_k(T)(a_1,\ldots, a_k))(a_{k+1})bc \\
&=&  a_0T(a_1,\ldots, a_k,a_{k+1}b)c - a_0T(a_1,\ldots,
a_k,a_{k+1})bc \\ & = & a_0T(a_1,\ldots, a_k,a_{k+1}b)c.
\end{eqnarray*}
However, $ a_0a_1=\cdots=a_k(a_{k+1}b)=(a_{k+1}b)c=0$, and $T$ is
$k+1$-\hl. Hence
$$a_0T(a_1,\ldots, a_k,a_{k+1}b)c=0.$$
Thus (1) holds, and so, by hypothesis, $S$ is a left multiplier.
Therefore $S(a)=S(1)a$ for all $a\in A$. However,
\begin{eqnarray*}
S(1) &=& [a_0\star \Lambda_k(T)(a_1,\ldots, a_k)\star a_{k+1}](1)
\\ &= & a_0(\Lambda_k(T)(a_1,\ldots, a_k))(a_{k+1}1)
- a_0(\Lambda_k(T)(a_1,\ldots, a_k))(a_{k+1})1 \\
&=& a_0T(a_1,\ldots, a_k,a_{k+1}) - a_0T(a_1,\ldots, a_k,a_{k+1})
\\ & = & 0.
\end{eqnarray*}
Thus $S=0$. Hence $\Lambda_k(T)$ is $k$-\hl. Let $q$ be the
natural quotient mapping from $B(A,X)$ into $B(A,X)/B_A(A,X)$,
where $B_A(A,X)$ is the space of left multipliers. Since
$\Lambda_k(T)$ is $k$-\hl\ and $q$ is an $A$-bimodule morphism with
the $\star$ actions, $q\circ \Lambda_k(T)$ is $k$-\hl. Moreover,
because of the assumption on $T$, $q\circ \Lambda_k(T)(a_1,\ldots,
a_k)=0$ if any one of $a_1,\ldots, a_k$ is 1. On the other hand,
for every $T\in B(A,X)$,
$$1\star T=T \ \ \text{and} \ T\star 1-T \in B_A(A,X).$$
Thus $B(A,X)/B_A(A,X)$ is a unital Banach $A$-bimodule. Therefore,
by the inductive hypothesis, $q\circ \Lambda_k(T)$ is a
$k$-cocycle. This means that for $a_1,\ldots, a_{k+1}\in A$,
$$\Delta^k(q\circ \Lambda_k(T))(a_1,\ldots, a_{k+1})=0.$$
Hence, from the equation $\Lambda_{k+1}\circ
\delta^{k+1}=\Delta^k\circ \Lambda_k$,
$$\Lambda_{k+1}(\delta^{k+1}(T))(a_1,\ldots, a_{k+1})=\Delta^k(\Lambda_k(T))(a_1,\ldots, a_{k+1})\in
B_A(A,X).$$ Thus, for every $a_{k+2}\in A$,
\begin{eqnarray*}
\delta^{k+1}(T)(a_1,\ldots, a_{k+1},a_{k+2}) &=&
[\Lambda_{k+1}(\delta^{k+1}(T))(a_1,\ldots, a_{k+1})](a_{k+2})
\\ &= & [\Lambda_{k+1}(\delta^{k+1}(T))(a_1,\ldots, a_{k+1})](1)a_{k+2}
\\ &= & \delta^{k+1}(T)(a_1,\ldots, a_{k+1},1)a_{k+2}.
\end{eqnarray*}
On the other hand, by the assumption on $T$,
$$a_1T(a_2,\ldots ,a_{k+1},1)+ \sum_{j=1}^k (-1)^j
T(a_1,\ldots ,a_ja_{j+1},\ldots ,a_{k+1},1)=0.$$ Also,
$$\delta^{k+1}(T)(a_1,\ldots,a_k, a_{k+1}1)-
\delta^{k+1}(T)(a_1,\ldots,a_k, a_{k+1})1=0.$$ Hence
$\delta^{k+1}(T)(a_1,\ldots, a_{k+1},1)=0$. Therefore
$\delta^{k+1}(T)=0$, and so $T\in B^{k+1}(A,X)$. This completes
the proof.
\end{proof}

We are now ready to state the main result of this section.
We recall that the unitization of A is $A^\sharp:=A \oplus \mathbb{C}$
with multiplication
 
$$(a, \lambda ) (b, \mu ) = ( ab + a\mu + b\lambda, \lambda \mu)  Ê ( a, b \in A,
\lambda, \mu \in \C ), $$ and norm
$$ \|( a, \lambda) \|= \| a \| +  |\lambda| \ \ \  ( a \in A,
\lambda \in \C ).$$ Thus $A^\sharp$ is a unital Banach algebra with
unit ($0,1$) which is denoted by 1 if there is no case of ambiguity. 

\begin{thm}\label{T:App Loc-n cocycles & ref}
Let $A$ be a Banach algebra such that $A^\sharp$  satisfies 
conditions $($i$)$ and $($ii$)$  of Proposition \ref{P:App. Loc. n
cocycle}. Then, for any Banach $A$-bimodule $X$ and $n\in \N$,
every bounded local $n$-cocycle $T$ from $A^{(n)}$
into $X$ is an $n$-cocycle. 
%In particular, $\mathcal{Z}^n(A,X)$ is
%reflexive and $\tref_a[Z^n(A,X)]\cap \mathcal{B}^n(A,X)=
%\mathcal{Z}^n(A,X)$.
\end{thm}

\begin{proof}
We can extend $X$ to a Banach $A^\sharp$-bimodule by defining
$1x=x1=x$. Let $\sigma \from L^n(A,X) \to L^n(A^\sharp ,X)$ be a
linear map defined by
$$\sigma(T)(a_1+\lambda_1,\ldots, a_n+\lambda_n)=T(a_1,\ldots,
a_n),$$ for $a_1,\ldots, a_n \in A$ and $\lambda_1,\ldots,
\lambda_n \in \C$. It is straightforward to check that $T\in
L^n(A,X)$ is an n-cocycle if and only if $\sigma(T)$ is an
n-cocycle. Now let $T\in B^n(A,X)$ be a bounded 
local n-cocycle, and let $(a_1+\lambda_1,\ldots, a_n+\lambda_n)\in
A^{\sharp (n)}$. By the assumption on $T$, for
$\tilde{a}=(a_1,\ldots, a_n)\in A^{(n)}$, there is an n-cocycle 
$T_{\tilde{a}}$ from $A^{(n)}$ into $X$ such that
$T(a_1,\ldots, a_n)=T_{\tilde{a}}(a_1,\ldots, a_n)$. Thus
\begin{eqnarray*}
\sigma(T)(a_1+\lambda_1,\ldots, a_n+\lambda_n) &=& T(a_1,\ldots, a_n) \\
&=& T_{\tilde{a}}(a_1,\ldots, a_n) \\&=&
\sigma(T_{\tilde{a}})(a_1+\lambda_1,\ldots,
a_n+\lambda_n).
\end{eqnarray*}
Hence $\sigma(T)$ is a bounded local
n-cocycle, and so it is a bounded n-hyperlocal map. 
Moreover, $$\sigma(T)(a_1,\ldots, a_n)=0$$ if any one of
$a_1,\ldots, a_n$ is 1. Thus, by Proposition \ref{P:App. Loc. n
cocycle}, $\sigma(T)$ is an n-cocycle. Therefore $T$ is an
n-cocycle.
\end{proof}

%\begin{rem}
%(i) It is not 

\section{Hyper-Tauberian algebras}\label{S:HL}

%In this section, we show that a hyper-Tauberian satisfies the assumption

Throughout this section, $A$ and $B$ are commutative regular
semisimple Banach algebras with the carrier spaces $\Phi_A$ and
$\Phi_B$, respectively. Let $I$ be a closed ideal in $A$. The {\it
hull} of $I$ is
$$\{t\in \Phi_A \mid a(t)=0\ \text{for all}\ a\in I \},$$ and it
is denoted by $h(I)$.

Let $X$ and $Y$ be Banach left (right) $A$-modules. For $x\in X$, the
annihilator Ann$_A$($x$) of $x$ is
$$\text{Ann}_A(x)=\{ a\in A \mid \ ax=0\ (xa=0)   \}. $$
Ann$_A$($x$) is clearly a closed ideal in $A$. The hull of
Ann$_A$($x$) is called the \emph{support} of $x$ (in $\Phi_A$),
denoted by $\supp_Ax$. We will write ``$\supp\,x$" instead of
``$\supp_Ax$" whenever there is no risk of ambiguity. By \cite[Lemma 2.1]{S1},
$t\notin \supp\,x$ if and only if there
is a compact neighborhood $V$ of $t$ in $\Phi_A$ such that, for
every element $a\in A$, if $\supp\,a \subseteq V$, then $ax=0 \
(xa= 0)$. In the case
$X=A$ where we regard $A$ as a Banach (left or right) $A$-module on
itself, the support of an element $a\in A$ coincides with the usual
definition of $\supp\,a$, namely cl$\{t\in \Phi_A \mid a(t)\neq 0
\}$.

 An operator
$T \from X \to Y$ is \emph{local}
with respect to the left (right) $A$-module action if $\supp\,T(x) \subseteq
\supp\,x$ for all $x\in X$. We recall from
\cite[Definition 4]{S3} that $A$ is a {\it \hT\ algebra} if every bounded local
operator from $A$ into $A^*$ is a multiplier. If $A$ is unital, then the definition
of \hT\ algebras coincides with the definition of (SD) algebras introduced
in \cite{Shul1}.

Let $\A$ and $\B$ be Banach algebras, and let $X$ be both a Banach
left $\A$-module and a Banach right $\B$-module such that for all
$a\in \A$, $b\in \B$ and $x\in X$, $a(xb)=(ax)b$. Then we write
$X\in \A-\tmod-\B$. If, in addition, $X$ is essential both as a
Banach left $\A$-module and Banach right $\B$-module, then we write
$X\in \text{ess.}\  \A-\tmod-\B$.

\begin{defn}\label{D:hlo}
Let $\A$ and $\B$ be Banach algebras, and let $X, Y\in
\A-\text{mod}-\B$. An operator $D \from X \to Y$ is {\it hyperlocal}
with respect to $\A$-mod-$\B$ actions if, for all $a\in \A$, $b\in \B$
and $x\in X$,
$$ax=xb=0 \ \ \  \text{implies} \ \ \ aD(x)b=0.$$
\end{defn}

The preceding definition was introduced in \cite{Sch} in order to
extend the concept of locality for operators in the
non-commutative setting (see also \cite{HL}). It is easily seen that, for commutative
C$^*$-algebras, this locality condition coincides with the usual
one. However, as it is shown in Remark \ref{R:lo-not hlo}, in
general the concept of being hyperlocal is weaker than the concept
of being local.

In the following proposition, we use the properties of hyper-Tauberian
algebras to characterize bounded hyperlocal operators that are defined
from essential modules over these algebras.

\begin{prop}\label{P:hT-hlo-Generalized der.}
Let $A$ and $B$ be \hT\ algebras. Then, for all $X,Z \in ess.\
A-mod-B$ and $Y\in ess.\ B-mod-A$,\\
$($i$)$ a bounded operator $D \from X \to Y^*$ is \hl\ if and only
if
$$D(axb)-aD(xb)-D(ax)b+aD(x)b=0$$ for all $a\in A$, $b\in B$ and
$x\in X$.\\
$($ii$)$ If $A$ and $B$ have bounded approximate identities, then
the result in $($i$)$ is also true for all bounded \hl\ operators
from $X$ into $Z$.
\end{prop}

\begin{proof}
(i) First assume that $Y=B\widehat{\otimes} A$, where the
$B$-mod-$A$ actions on $B \widehat{\otimes} A$ are specified by
$$ d (b \otimes a)=db \otimes a \ \ , \ \ (b \otimes a)
c=b \otimes ac \ \ (a,c\in A, b,d\in B).$$ Let $D \from X \to
(B\widehat{\otimes} A)^*$ be a bounded \hl\ operator, and let $x\in
X$ and $a\in A$. Define the bounded operator $\widetilde{D} \from B
\to (B\widehat{\otimes} A)^*$ by
$$\widetilde{D}(b)=D(axb)-aD(xb) \ \ (b\in B).$$
We claim that $D$ is local with respect to right $B$-module action.
Let $b\in B$ and $t\notin \supp_B\, b$. There is a compact
neighborhood $V$ of $t$ (in $\Phi_B$) such that $V\cap \supp_B\,
b=\emptyset$. Let $c\in B$ with $\supp_B\,c \subseteq V$. By the
regularity of $B$, there is $e\in B$ such that $e=1$ on $V$ and
$e=0$ on $\supp_B\,b$. So $$ec=c \ \text{and}\ eb=0. \eqno{(1)}$$
Put
$$K_0(V)=\overline{\text{span}}\{ n\otimes m \mid m\in A\ ,n\in B\
\text{and}\ n=0 \ \text{on}\ \Phi_B \setminus V \}.$$ Since $e=1$ on
$V$, for all $\theta\in (B\widehat{\otimes} A)^*$,
$$\theta e-\theta=0\ \text{on}\ K_0(V).\eqno{(2)}$$ Let $z\in X$,
and define the bounded operator $T \from A \to (B\widehat{\otimes}
A)^*/ K_0(V)^\perp$ by
$$T(u)=D(uzb)+K_0(V)^\perp \ \ (u\in A).$$
Let $h\in A$ such that $hu=0$. Then, from (1), $huzb=0=uzbe$. Since
$D$ is \hl, $hD(uzb)e=0$. Hence, from (2),
\begin{eqnarray*}
hT(u) &=& hD(uzb)+K_0(V)^\perp \\ &=& hD(uzb)e+K_0(V)^\perp \\ & = &
0.
\end{eqnarray*}
In particular, $T$ is local with respect to left $A$-module action.
Since $$(B\widehat{\otimes} A)^*/ K_0(V)^\perp\cong K_0(V)^*,$$ and
$K_0(V)$ is an essential Banach right $A$-module, from \cite[Proposition
3]{S3}, it follows that $T$ is a right multiplier.
Therefore $T(uv)=uT(v)$ for all $u, v\in A$. Hence, if we put $u=a$,
then $D(avzb)-aD(vzb)\in K_0(V)^\perp$. Thus, from essentiality of
$X$, we have
$$\widetilde{D}(b)=D(axb)-aD(xb)\in K_0(V)^\perp.$$
Therefore $\widetilde{D}(b)c=0$, since $\supp_B\, c\in V$. This
means that $t\notin \supp_B\, \widetilde{D}(b)$, and so
$\widetilde{D}$ is a bounded local operator. Hence, from \cite[Proposition 3]{S3},
$\widetilde{D}$ is a left multiplier. Thus
$\widetilde{D}(bd)=\widetilde{D}(b)d$ for all $b,d\in B$. Therefore
$$D(axbd)-aD(xbd)=D(axb)d-aD(xb)d.$$
The final result follows from the essentiality of $X$.\\
Now consider the general case. Let $y\in Y$ and define $S_y \from
Y^* \to (B\widehat{\otimes} A)^*$ by
$$\la S_y(y^*)\ , \ b\otimes a \ra= \la y^*\ , \ bya \ra \ \ \ \ (a\in A \ , \ b\in B \ , \ y^*\in Y^*).$$It
is easy to see that $S_y$ is both a bounded left $A$-module morphism
and a bounded right $B$-module morphism, and so $S_x\circ D$ is a
bounded hyperlocal operator from $X$ into $(B\widehat{\otimes}
A)^*$. Thus, for all $a\in A$, $b\in B$, $x\in X$ and $y\in Y$,
$$ S_y[D(axb)-aD(xb)-D(ax)b+aD(x)b]=0.$$
Hence, for all $c\in A$ and $d\in B$,
$$ \la D(axb)-aD(xb)-D(ax)b+aD(x)b\ , \ dyc \ra =0.$$
The final results follows from the essentiality of $Y$,\\
(ii)  Let $\{e_\alpha \}_{\alpha \in \Lambda}$ and $\{f_\beta
\}_{\beta \in \Omega}$ be bounded approximate identities for $A$ and
$B$, respectively. Similar to the argument made in (i) (by replacing
$Z$ with $Z^{**}$), we can show that
$$c[D(axb)-aD(xb)-D(ax)b+aD(x)b]d=0 \eqno{(3)}$$ for all $a,c \in A$, $b,d \in B$
and $x\in X$. On the other hand, since $A$ and $B$ have bounded
approximate identities, by Cohen's factorization theorem
\cite[Theorem 11.10]{BD}, there are $e\in A$, $f\in B$ and $z\in Z$
such that $$D(axb)-aD(xb)-D(ax)b+aD(x)b=ezf.$$ So we have the final
result if we put $c=e_{\alpha}$ and $d=f_{\beta}$ in (3), and let
$\alpha, \beta \to \infty$.
\end{proof}

\begin{thm}\label{T:hT-ref-n cocycles}
Let $A$ be a \hT\ algebra, and let $X$ be a Banach $A$-bimodule.
Then, for $n\in \N$, every bounded local $n$-cocycle
$T$ from $A^{(n)}$ into $X$ is an $n$-cocycle.
 %In particular, $\mathcal{Z}^n(A,X)$ is reflexive.
\end{thm}

\begin{proof}
Let $A^\sharp$ be the unitalization of $A$. By \cite[Corollary 10]
{S3}, $A^\sharp$ is \hT. Therefore, by
\cite[Proposition 3]{S3} and Proposition
\ref{P:hT-hlo-Generalized der.}, $A^\sharp$ satisfies the
conditions (i) and (ii) of Proposition \ref{P:App. Loc. n
cocycle}. Hence the result follows from Theorem \ref{T:App Loc-n
cocycles & ref}.
\end{proof}

\begin{rem}\label{R:lo-not hlo}
Let $\mathbb{T}$ be the unit circle, and let
$A:=A(\mathbb{T})$ be the Fourier algebra on $\mathbb{T}$. It is
shown in \cite[Remark 24(ii)]{S3} that there is a closed ideal
$I$ in $A$ such that $I$ is weakly amenable but $I$ is not \hT.
Hence there are bounded local operators from $I$ into $I^*$ which
are not multipliers. However, this is not the case if we consider
bounded hyperlocal operators. To see this, let $D \from I \to I^*$
be a bounded hyperlocal operator.  First, we show that $D$ is
hyperlocal with respect to $A$-bimodule actions. Let $a,b\in A$
and $c\in I$ such that $ac=cb=0$. Take $e,f\in I$. Then $ea,bf\in
I$ and $(ea)c=c(bf)=0$. Thus $eaD(c)bf=0$. Hence $aD(c)b=0$ on
$I^3$. However, $I$ is weakly amenable, and so, by \cite[Theorem
2.8.69(ii)]{D}, $I^2$ is dense in $I$. Hence $aD(c)b=0$. Therefore
$D$ is hyperlocal with respect to $A$-bimodule actions. On the
other hand, $A$ is a \hT\ algebra \cite[Proposition 18]{S3}.
Hence, from Proposition \ref{P:hT-hlo-Generalized der.}, for all
$a,b\in A$ and $c\in I$,
$$D(abc)-aD(bc)-D(ab)c+aD(b)c=0. \eqno{(1)}$$
Define the bounded operator $\mathcal{D} \from I \to
\mathcal{B}_I(I,I^*)$ by
$$ \mathcal{D}(a)(b)=D(ab)-aD(b) \ \ (a,b\in I). $$
From (1), it is easy to verify that $\mathcal{D}$ is well-defined.
Moreover, upon setting
$$\la a\cdot S \ , \ b \ra =\la S\cdot a \ , \ b\ra =\la S \ , \ ab \ra,$$ the space
$\mathcal{B}_I(I,I^*)$ becomes a symmetric Banach $I$-module and
$\mathcal{D}$ becomes a bounded derivation from $I$ into
$\mathcal{B}_I(I,I^*)$. Hence $\mathcal{D}=0$ since $I$ is weakly
amenable. Thus $D$ is a multiplier.
\end{rem}

\section{C$^*$-algebras}\label{S:C alg}

It follows from the works of B. E. Johnson that $C_0(\mathbb{R})$ is
a hyper-Tauberian algebra \cite[Proposition 3.1]{J}. One the other hand,
Shulman showed that every unital commutative C$^*$-algebra is hyper-Tauberian \cite{Shul1}.
We extend these results by showing that $C_0(\Omega)$ is
hyper-Tauberian for every locally compact topological space $\Omega$.
For the sake of completeness, we first prove it for the
case when $\Omega$ is compact.

\begin{thm}\label{T:hT-Comm. C*-al}
Let $\Omega$ be a locally compact topological space. Then
$C_0(\Omega)$ is a \hT\ algebra.
\end{thm}

\begin{proof}
First consider the case when $\Omega$ is compact. Let $T \from
C(\Omega) \to C(\Omega)^*$ be a bounded local operator. First we
show that $T$ satisfies the following condition:
$$ ab=0\ \ \text{implies}\ \ aT(b)=0. \eqno{(\star)}$$
Let $a,b\in C(\Omega)$ with $ab=0$. So if we put $E=\supp\,b$, then
$a=0$ on $E$. Since $E$ is a closed subset of $\Omega$, $E$ is a
set of synthesis (see \cite[Definition 4.1.12 and Theorem 4.2.1]{D}). Thus there is a sequence
$\{a_n\}$ in $C(\Omega)$ such that, for each $n$, $\supp\, a_n$ is 
compact and disjoint from $E$, and $a_n \to a$ as $n \to \infty$. On
the other hand, since $T$ is local and $\supp\,a_n$ is disjoint
from $E$,
\begin{eqnarray*}
\supp\,a_n \cap \supp\,T(b) &\subseteq & \supp\,a_n \cap \supp\,b
\\ &=& \supp\,a_n\cap  E \\ & = & \emptyset.
\end{eqnarray*}
Therefore, since $\supp\,a_n$ is compact, $a_nT(b)=0$. Hence, by
letting $n \to \infty$, we have $aT(b)=0$. This proves ($\star$).
Now let $a\in C(\Omega)$ be a self-adjoint element, and let $A(a)$
be the $C^*$-subalgebra of $C(\Omega)$ generated by $\{a,1\}$. It
is well-known that there is a compact subset $K$ of $\mathbb{R}$
such that $A(a)$ is isometrically isomorphic to $C(K)$. In
particular, $C(\Omega)$ is an essential and symmetric Banach
$C(K)$-module. Let $d\in C(\Omega)$ and $c\in C(K)$ with $cd=0$.
Then, since $c\in A$ and $T$ satisfies condition ($\star$),
$cT(d)=0$. Hence Ann$_{C(K)}\,d \subseteq
\text{Ann}_{C(K)}\,T(d)$, and so $\supp_{C(K)}\,T(d)\subseteq
\supp_{C(K)}\,d$. Therefore $T$ is local with respect to
$C(K)$-module actions. On the other hand,
the restriction map $f \mapsto f|_K$ is a
bounded algebra homomorphism from $C_0(\mathbb{R})$ onto $C(K)$.
Hence, from
\cite[Theorem 12]{S3}, $C(K)$ is \hT, and so, from
\cite[Proposition 3]{S3}, $T$ is a $C(K)$-module morphism.
Hence, for each $b\in C(\Omega)$, $T(ab)=aT(b)$. The final result
follows since $C(\Omega)$ is the linear span of its self-adjoint
elements.\\
We now consider the general case. Let $\Omega$ be a locally
compact topological space, and let $\Omega\cup \{\infty \}$ be its one-point
compactification. Then, from the first case, $C(\Omega\cup
\{\infty \})$ is \hT. On the other hand,
$$C_0(\Omega)=\{a\in C(\Omega\cup \{\infty \}) \mid
a(\infty)=0 \}$$ and $\{\infty \}$ is a set of
synthesis for $C(\Omega\cup \{\infty \})$. Thus, from \cite[Theorem
7(ii)]{S3}, $C_0(\Omega)$ is
\hT.
\end{proof}

We are now ready to obtain our results for C$^*$-algebras.
We start with the following critical theorem which characterizes 
bounded hyperlocal operators defined over essential modules of 
a C$^*$-algebra. This was partially obtained in \cite[Theorem 2.2]{Sch}
and \cite[Theorem 2.17]{HL}.

\begin{thm}\label{T:hlo-C* alg}
Let $A$ be a $C^*$-algebra, let $X$ be an essential Banach
$A$-bimodule, and let $Y$ be an essential or the dual of an
essential Banach $A$-bimodule. Then a bounded operator $D \from X
\to Y$ is hyperlocal if and only if
$$D(axb)-aD(xb)-D(ax)b+aD(x)b=0$$ for all $a,b\in A$ and $x\in X$.
\end{thm}

\begin{proof}
First assume that $Y=(A\widehat{\otimes} A)^*$. Let $D \from X \to
(A\widehat{\otimes} A)^*$ be a bounded hyperlocal operator, and
let $A^\sharp$ be the unitalization of $A$ \cite[Definition
3.2.1]{D}. We show that $D$ is hyperlocal with respect to
$A^\sharp$-module actions. Let $u,v\in A^\sharp$ and $x\in X$ such
that $ux=xv=0$. So, for all $a,b\in A$, $(au)x=x(vb)=0$. Thus
$auD(x)vb=0$. Hence $uD(x)v=0$ on $A^2 \otimes A^2$ which is dense
in $A\widehat{\otimes} A$. So $uD(x)v=0$. Now let $c$ and $d$ be
self-adjoint elements in $A$, and let $A(c)$ and $A(d)$ be the
commutative $C^*$-subalgebras of $A^\sharp$ generated by $\{c,1\}$ and
$\{d,1\}$, respectively. Clearly $D \from X \to (A\widehat{\otimes} A)^*$ is \hla\ with
respect to $A(c)$-mod-$A(d)$ actions. Thus, from Theorem \ref{T:hT-Comm. C*-al}, for
every $x\in X$,
$$D(cxd)-cD(xd)-D(cx)d+cD(x)d=0.$$
The final result follows since $A$ is the linear span of its
self-adjoint elements.
The general case follows from a similar argument made in the proof of
Proposition~\ref{P:hT-hlo-Generalized der.}.
\end{proof}

\begin{rem}\label{R:C* alg-left hlo}
In the preceding theorem, if we replace the locality condition
that we used in the definition of a \hla\ operator with the
following condition
$$ ax=0 \ \ \text{implies}\ \ aD(x)=0,$$
then, by a similar argument and using \cite[Proposition 3]{S3}
instead of Proposition \ref{P:hT-hlo-Generalized der.}, we can show that $D$ is a left $A$-module morphism. We can
also have a similar result regrading bounded right $A$-module
morphisms.
\end{rem}

Let $A$ be a $C^{\ast}$-algebra which is not unital. We can see
that, in general, our unitization, $A^\sharp=A\oplus^1 \C$, is not
a $C^{\ast}$-algebra (as the norm dose not satisfy the correct
condition). However, there is an equivalent norm on $A^\sharp$
that turns it into a $C^*$-algebra (see \cite[Definition
3.2.1]{D}). Thus we can state the our main result:

\begin{thm}\label{T:C*-ref-n cocyclyes}
Let $A$ be a $C^*$-algebra, and let $X$ be a Banach $A$-bimodule.
Then, for $n\in \N$, every bounded local $n$-cocycle
$T$ from $A^{(n)}$ into $X$ is an $n$-cocycle. 
%In particular,
%$\mathcal{Z}^n(A,X)$ is reflexive and $\tref_a[Z^n(A,X)]\cap \mathcal{B}^n(A,X)=
%\mathcal{Z}^n(A,X)$.
\end{thm}

\begin{proof}
The result follows from Theorem \ref{T:hlo-C* alg}, Remark \ref{R:C* alg-left hlo}
and Theorem \ref{T:App Loc-n cocycles & ref}.
\end{proof}

\section{Hyperlocal operators and Amenable C$^*$-algebras}

In this final section, we present a characterization of
cohomological properties of C$^*$-algebras, i.e. amenability and weak amenability,
with respect to hyperlocal operators.

\begin{thm}\label{T:hlo-C* alg-dual module}
Let $A$ be a $C^*$-algebra, and let $X$ be an essential Banach
$A$-bimodule. Then a bounded operator $D \from A \to X^*$ is
hyperlocal if and only if there is a derivation $\mathcal{D}$ and a
right multiplier $T$ from $A$ into $X^*$ such that
$D=\mathcal{D}+T$. In particular, $D$ is a derivation if and only if
$weak^*-\dss\lim_{\alpha\to\infty} D(e_{\alpha})=0$ for a bounded
approximate identity $\{e_{\alpha}\}_{\alpha \in \Lambda}$ in $A$.
\end{thm}

\begin{proof}
It is easy to see that all derivations and multipliers are
hyperlocal. On the other hand, let $D \from A \to X^*$ be a bounded
hyperlocal operator. By Theorem \ref{T:hlo-C* alg}, for
all $a,b,c\in A$,
$$D(acb)-D(ac)b-aD(cb)+aD(c)b=0.$$
By putting $c=e_{\alpha}$ and letting $\alpha \to \infty$ we obtain
$$D(ab)-D(a)b-aD(b)+\dss\lim_{\alpha\to\infty} aD(e_{\alpha})b=0. $$
Since $\{D(e_{\alpha})\}$ is bounded, there is $x^*\in X^*$ and a
subnet $\{D(e_{\alpha_i})\}$ such that $D(e_{\alpha_i}) \to x^*$ in
the weak$^*$ topology. So $D(ab)-D(a)b-aD(b)+ax^*b=0$. Define $T
\from A \to X^*$ by $T(a)=ax^*$ and put $\mathcal{D}=D-T$. It is
straightforward to check that $T$ is a right multiplier and
$\mathcal{D}$ is a derivation. Finally, $D$ is a derivation if and
only if $T$ is zero. However, it is easy to verify that $T$ is zero
if and only if $\text{weak}^*-\dss\lim_{\alpha\to\infty}
D(e_{\alpha})=0$.
\end{proof}

We recall that a Banach algebra $A$ is amenable if  for
any Banach $A$-bimodule $X$, every bounded derivation 
$D \from A \to X^*$ is inner. 

\begin{cor}\label{C:hlo-C* al-amenable}
Let $A$ be a $C^*$-algebra. Then $A$ is amenable if and only if for
any essential Banach $A$-bimodule $X$ and every bounded hyperlocal
operator $D \from A \to X^*$, there are $x^*, y^* \in X^*$ such that
$D(a)=ax^*-y^*a$ $(a\in A)$.
\end{cor}

\begin{proof}
Let $A$ be amenable, let $X$ be an essential Banach $A$-bimodule,
and let $D \from A \to X^*$ be a bounded \hla\ operator. By
Theorem \ref{T:hlo-C* alg-dual module}, there is a derivation $\mathcal{D}$
and a right multiplier $T$ from $A$ into $X^*$ such that
$D=\mathcal{D}+T$. Since $A$ is amenable, there are $y^*$ and $z^*$
in $X^*$ such that $\mathcal{D}(a)=ay^*-y^*a$ and $T(a)=az^*$ for
all $a\in A$. Thus $D(a)=a(y^*+z^*)-y^*a$. The converse follows
immediately from Theorem \ref{T:hlo-C* alg-dual module} and \cite[Corollary
2.9.27]{D}.
\end{proof}

\begin{cor}
Let $A$ be a $C^*$-algebra. Then, for every bounded hyperlocal
operator $D \from A \to A^*$, there are $x^*, y^* \in A^*$ such that
$D(a)=ax^*-y^*a$ $(a\in A)$.
\end{cor}

\begin{proof}
The result follows from the similar argument to the one made in
the proof of the preceding corollary together with the fact that every
$C^*$-algebra is weakly amenable \cite[Theorem 5.6.77]{D}.
\end{proof}


\begin{thebibliography}{99}
\bibitem{BD} F. F. Bonsall and J. Duncan, \emph{Complete Normed Algebras},
New York, Springer-Verlag 1973.
\bibitem{D} H. G. Dales, {\it Banach algebras and automatic
continuity}, New York, Oxford University Press, 2000.
\bibitem{HL} D. Hadwin and J. Li, {\it Local derivations and local
automorphism}, J. Math. Anal. Appl. {\bf 290} (2004) 702-714.
\bibitem{J} B. E. Johnson, \emph{Local derivations on $C^*$-algebras are derivations},
Trans. Amer. Math. Soc, \textbf{353 }(2000), 313-325.
\bibitem{K} R. V. Kadison,\emph{ Local derivations}, J. Algebra \textbf{130}
(1990), 494-509.
\bibitem{L} David R. Larson, {\it Reflexivity, algebraic
reflexivity and linear interpolation}, Amer. J. Math. {\bf 110}
(1988), 283-299.
\bibitem{S1} E. Samei, {\it Bounded and completely bounded local derivations from
certain commutative semisimple Banach algebras}, Proc. Amer. Math.
Soc, {\bf 133} (2005), 229-238.
\bibitem{S2} E. Samei, {\it Approximately local derivations},
J. London Math. Soc. (2) {\bf 71} (2005), no. 3, 759--778.
\bibitem{S3} E. Samei, {\it Hyper-Tauberian algebras and weak amenability of
Fig\`{a}-Talamanca-Herz algebras}, J. Func. Anal.  {\bf 231}  (2006), no. 1, 195-220.
\bibitem{Sch} J. Schweizer, {\it An analogue of Peetr's theorem in
non-commutative topology}, Quart. J. Math. {\bf 52} (2001)
499-506.
\bibitem{Shul1} V. S. Shulman, {\it Spectral synthesis and the Fuglede-Putnam-Rosenblum
theorem.} (Russian) Teor. Funktsi\u\i Funktsional. Anal. i Prilozhen. No. 54 (1990), 25--36;
translation in J. Soviet Math. 58 (1992), no. 4, 312--318.
%\bibitem{T} J. Tomiyama, \emph{Tensor products of commutative Banach
%algebras}, T\^{o}hoku Math. J. (2)\textbf{ 12} (1960),
%143-154.
\end{thebibliography}
\end{document}